\input amstex
\input Amstex-document.sty

\pageno 173

\def\shorttitle{Eisenstein Series and Arithmetic Geometry}
\def\shortauthor{S. S. Kudla}

\define\C{{\Bbb C}}

\define\R{{\Bbb R}}
\define\Q{{\Bbb Q}}
\define\Z{{\Bbb Z}}

\redefine\H{\frak H}

\define\a{\alpha}

\define\e{\epsilon}

\redefine\l{\lambda}
\redefine\o{\omega}

\redefine\P{\Phi}
\predefine\Sec{\S}
\redefine\L{\Lambda}

\define\back{\backslash}
\define\lra{\longrightarrow}
\define\scr{\scriptstyle}
\define\liminv#1{\underset{\underset{#1}\to\leftarrow}\to\lim}
\define\limdir#1{\underset{\underset{#1}\to\rightarrow}\to\lim}

\define\nass{\noalign{\smallskip}}

\define\CH{\widehat{CH}}


\redefine\ord{\text{\rm ord}}
\define\Ei{\text{\rm Ei}}
\redefine\O{\Omega}
\predefine\oldvol{\vol}
\redefine\vol{\text{\rm vol}}

\redefine\div{\text{\rm div}}

\define\Spec{\text{\rm Spec}\,}

\define\Sym{\text{\rm Sym}}
\define\tr{\text{\rm tr}}

\font\cute=cmitt10 at 12pt
\font\smallcute=cmitt10 at 9pt
\define\kay{{\text{\cute k}}}
\define\smallkay{{\text{\smallcute k}}}
\define\OK{\Cal O_{\smallkay}}

\define\End{\text{\rm End}}

\define\SL{\text{\rm SL}}

\define\Pic{\text{\rm Pic}}
\define\Pich{\widehat{\Pic}}
\define\degh{\widehat{\deg}\ }
\define\divh{\widehat{\div}\ }
\define\GSpin{\text{\rm GSpin}}

\define\hfb{\hfill\break}

\define\Zh{\widehat{\Cal Z}}

\define\borchduke{1}
\define\bost{2}
\define\boutotcarayol{3}
\define\brkuehn{4}
\define\bkk{5}
\define\gsihes{6}
\define\grosskeating{7}
\define\kitaoka{8}
\define\duke{9}
\define\annals{10}
\define\kbourb{11}
\define\kudlamsri{12}
\define\kmihes{13}
\define\krHB{14}
\define\krinvent{15}
\define\krsiegel{16}
\define\tiny{17}
\define\kryII{18}
\define\kuehn{19}
\define\mcgraw{20}
\define\shimura{21}
\define\yangden{22}
\define\yangiccm{23}
\define\yangmsri{24}
\define\zagier{25}

\topmatter %
\title\nofrills{\boldHuge Derivatives of Eisenstein Series and Arithmetic Geometry*}
\endtitle

\author \Large Stephen S. Kudla$^\dag$\endauthor

\thanks *Partially
supported by NSF grant DMS-9970506 and by a Max-Planck Research Prize from the Max-Planck Society and Alexander
von Humboldt Stiftung. \endthanks

\thanks $^\dag$ Mathematics Department, University of Maryland, College Park, MD 20742, USA. E-mail:
ssk\@math.umd.edu
\endthanks

\abstract\nofrills \centerline{\boldnormal Abstract}

\vskip 4.5mm

{\ninepoint We describe connections between the Fourier coefficients of derivatives of Eisenstein series and
invariants from the arithmetic geometry of the Shimura varieties $M$ associated to rational quadratic forms
$(V,Q)$ of signature $(n,2)$. In the case $n=1$, we define generating series $\hat\phi_1(\tau)$ for $1$-cycles
(resp. $\hat\phi_2(\tau)$ for $0$-cycles) on the arithmetic surface $\Cal M$  associated to a Shimura curve over
$\Bbb Q$. These series are related to the second term in the Laurent expansion of an Eisenstein series of weight
$\frac32$ and genus $1$ (resp. genus $2$) at the Siegel--Weil point, and these relations can be seen as examples
of an `arithmetic' Siegel--Weil formula. Some partial results and conjectures for higher dimensional cases are
also discussed.

\vskip 4.5mm


}
\endabstract
\endtopmatter

\document

\baselineskip 4.5mm \parindent 8mm

\specialhead \noindent \boldLARGE 1. Introduction \endspecialhead

In this report, we will survey results about generating functions for arithmetic cycles on Shimura varieties
defined by rational quadratic forms of signature $(n,2)$. For small values of $n$, these Shimura varieties are of
PEL type, i.e., can be identified with moduli spaces for abelian varieties equipped with polarization,
endomorphisms, and level structure. By analogy with CM or Heegner points on modular curves, cycles are defined by
imposing additional endomorphisms. Relations between the heights or arithmetic degrees of such cycles and the
Fourier coefficients of {\it derivatives} of Siegel Eisenstein series are proved in \cite{\annals} and in
subsequent joint work with Rapoport, \cite{\krHB}, \cite{\krinvent},  \cite{\krsiegel}, and with Rapoport and Yang
\cite{\tiny}, \cite{\kryII}. These relations may be viewed as an arithmetic version of the classical Siegel-Weil
formula, which identifies the Fourier coefficients of {\it values} of Siegel Eisenstein series with representation
numbers of quadratic forms. The most complete  example is that of anisotropic ternary quadratic forms ($n=1$), so
that the cycles are curves and $0$-cycles on the arithmetic surfaces associated to Shimura curves. Other surveys
of the material discussed here can be found in \cite{\kbourb} and \cite{\kudlamsri}.

\specialhead \noindent \boldLARGE 2. Shimura curves \endspecialhead

Let $B$ be an indefinite quaternion algebra over $\Q$, and
let $D(B)$ be the product of the primes $p$ for which
$B_p = B\otimes_\Q\Q_p$ is a division algebra. The rational vector space
$$V= \{\ x\in B\mid \tr(x)=0\ \}$$
with quadratic form given by $Q(x) = -x^2 =\nu(x)$,
where
$\tr(x)$ (resp. $\nu(x)$) is the reduced trace (resp. norm)
of $x$, has signature $(1,2)$. The
action of $B^\times$ on $V$ by conjugation
gives an isomorphism $G= \GSpin(V)\simeq B^\times$. Let
$$D = \{\ w\in V(\C)\mid (w,w)=0, \ (w,\bar w)<0\ \}/\C^\times
\simeq \Bbb P^1(\C)\setminus \Bbb P^1(\R)$$ be the associated symmetric space. Let $O_B$ be a maximal order in $B$
and let $\Gamma = O_B^\times$ be its unit group. The quotient $M(\C) = \Gamma\back D$ is the set of complex points
of the Shimura curve $M$ (resp. modular curve, if $D(B)=1$) determined by $B$. This space should be viewed as an
orbifold $[\Gamma\back D]$. For a more careful discussion of this and of the stack aspect, which we handle loosely
here, see \cite{\kryII}. The curve $M$ has a canonical model over $\Q$. From now on, we assume that $D(B)>1$, so
that $M$ is projective. Drinfeld's model $\Cal M$ for $M$ over $\Spec(\Z)$ is obtained as the moduli stack for
abelian schemes $(A,\iota)$ with an action $\iota:O_B\hookrightarrow \End(A)$ satisfying the `special' condition,
\cite{\boutotcarayol}. It is proper of relative dimension $1$ over $\Spec(\Z)$, with semi-stable reduction at all
primes and is smooth at all primes $p$ at which $B$ splits, i.e., for $p\nmid D(B)$. We view $\Cal M$ as an
arithmetic surface in the sense of Arakelov theory and consider its arithmetic Chow groups with real coefficients
$\CH^r(\Cal M) = \CH^r_{\R}(\Cal M)$, as defined in \cite{\bost}. Recall that these groups are generated by pairs
$(\Cal Z,g)$, where $\Cal Z$ is an $\R$-linear combination of divisors on $\Cal M$ and $g$ is a Green function for
$\Cal Z$, with relations given by $\R$-linear combinations of elements $\divh(f) = (\div(f),-\log|f|^2)$ where
$f\in \Q(\Cal M)^\times$ is a nonzero rational function on $\Cal M$. These real vector spaces come equipped with a
geometric degree map $\deg_\Q:\CH^1(\Cal M) \rightarrow CH^1(\Cal M_\Q) \overset{\deg}\to{\rightarrow} \R$, where
$\Cal M_\Q$ is the generic fiber of $\Cal M$, an arithmetic degree map $\degh:\CH^2(\Cal M)\rightarrow \R$,  and
the Gillet-Soul\'e height pairing, \cite{\bost},
$$\langle \ ,\ \rangle :\CH^1(\Cal M)\times \CH^1(\Cal M) \lra \R.$$
Let $\Cal A$ be the universal abelian scheme over $\Cal M$. Then
the Hodge line bundle
$\o = \e^*(\O^2_{\Cal A/\Cal M})$
determined by $\Cal A$ has
a natural metric, normalized as in \cite{\kryII}, section 3, and defines an
element $\hat\o\in \Pich(\Cal M)$, the group of metrized line bundles on $\Cal M$.
We also write $\hat\o$ for the image
of this class in $\CH^1(\Cal M)$ under the natural map, which sends
a metrized line bundle $\hat\Cal L = (\Cal L, ||\ ||)\in \Pich(\Cal M)$
to the class of $(\div(s),-\log||s||^2)$, for any nonzero section $s$ of $\Cal L$.

Arithmetic cycles in $\Cal M$ are defined by imposing additional endomorphisms of the following type.
\proclaim{Definition 1} {\rm (\cite{\annals})} The space of special endomorphisms $V(A,\iota)$ of  an abelian
scheme $(A,\iota)$, as above, is
$$V(A,\iota) = \{\ x\in \End(A)\mid x\circ\iota(b) = \iota(b)\circ x,\
\forall\, b\in O_B,\
\text{\rm and}\  \tr(x)=0 \ \},$$
with $\Z$-valued quadratic form given by $-x^2=Q(x)\,\roman{id}_A$.
\endproclaim

\specialhead \noindent \boldlarge 2.1. Divisors\endspecialhead To obtain divisors on $\Cal M$, we impose a single
special endomorphism.  For a positive integer $t$, let $\Cal Z(t)$ be the divisor on $\Cal M$ determined by the
moduli stack of triples $(A,\iota,x)$ where $(A,\iota)$ is as before and where $x\in V(A,\iota)$ is a special
endomorphism with $Q(x)=t$. Note that, for example, the complex points $\Cal Z(t)(\C)$ of $\Cal Z(t)$ correspond
to abelian surfaces $(A,\iota)$ over $\C$ with an `extra' action of the order $\Z[\sqrt{-t}]$ in the imaginary
quadratic field $\Q(\sqrt{-t})$, i.e., to CM points on the Shimura curve $M(\C)$. On the other hand, the cycles
$\Cal Z(t)$ can have vertical components in the fibers of bad reduction $\Cal M_p$ for $p\mid D(B)$. More
precisely, in joint work with M. Rapoport we show: \proclaim{Proposition 1} {\rm (\cite{\krinvent})} For $p\mid
D(B)$, $\Cal Z(t)$ contains components of the  fiber of bad reduction $\Cal M_p$ if and only if $\ord_p(t)\ge2$
and no prime $\ell\mid D(B)$, $\ell\ne p$, splits in $\kay_t:=\Q(\sqrt{-t})$.
\endproclaim
The precise structure of the vertical part of $\Cal Z(t)$ is determined in
\cite{\krinvent} using the Drinfeld-Cherednik $p$-adic uniformization of $\Cal M_p$.
For example, for $p\mid D(B)$, the multiplicities of the vertical components
in the fiber $\Cal M_p$ of the cycle $\Cal Z(p^{2r}t)$ grow with $r$, while the
horizontal part of this cycle remains unchanged.

To obtain classes in $\CH^1(\Cal M)$,
we construct Green functions by the procedure introduced in \cite{\annals}.
Let $L=O_B\cap V$.  For $t\in\Z_{>0}$
and $v\in \R_{>0}$,   define
a function $\Xi(t,v)$ on $M(\C)$ by
$$\Xi(t,v)(z) = \sum_{x\in L(t)} \beta_1(2\pi v R(x,z)),$$
where $L(t)=\{x\in L\mid Q(x)=t\}$, and, for $z\in D$ with preimage
$w\in V(\C)$,
$R(x,z) = |(x,w)|^2\,|(w,\bar w)|^{-1}$.
Here
$$\beta_1(r) = \int_1^\infty e^{-r u} \,u^{-1}\,du = -\Ei(-r)$$
is the exponential integral. Recall that this function has a log singularity
as $r$ goes to zero and decays exponentially as $r$ goes to infinity.
In fact, as shown in \cite{\annals}, section 11, for any  $x\in V(\R)$
with $Q(x)\ne0$,
the function
$$\xi(x,z) := \beta_1(2\pi R(x,z))$$
can be viewed as a Green function on $D$ for the divisor
$D_{x} :=\{z\in D\mid (x,z)=0\}.$
A simple calculation, \cite{\annals}, shows that, for $t>0$, $\Xi(t,v)$
is a Green function of logarithmic
type for the cycle $\Cal Z(t)$, while, for $t<0$, $\Xi(t,v)$ is a
smooth function on $M(\C)$.
\proclaim{Definition 2}{\rm(}i\,{\rm)} For $t\in \Z$ and $v>0$, the class $\Zh(t,v)\in
\CH^1(\Cal M)$  is defined by:
$$\Zh(t,v) = \cases (\Cal Z(t),\Xi(t,v))&\text{ if $t>0$,}\\
\nass
-\hat\o + (0,\bold c-\log(v))&\text{ if $t=0$, }\\
\nass
(0,\Xi(t,v))&\text{ if $t<0$.}
\endcases
$$
Here $\hat\o$ is the metrized Hodge line bundle, as above, and the real constant $\bold c$
is given by
$$\frac12\,\deg_\Q(\hat\o)\cdot \bold c
= \langle \hat\o,\hat\o\rangle - \zeta_{D(B)}(-1)\,\left[2\,\frac{\zeta'(-1)}{\zeta(-1)}
+ 1-\log(4\pi) -\gamma -\sum_{p\mid D(B)} \frac{p\log(p)}{p-1}\,\right],$$
where $\zeta_{D(B)}(s) = \zeta(s)\,\prod_{p\mid D(B)} (1-p^{-s})$ and $\gamma$ is Euler's
constant. \hfb
{\rm(}ii\,{\rm)} For $\tau = u+iv\in \H$ and $q= e(\tau) = e^{2\pi i\tau}$, the
`arithmetic theta function'
$\hat\phi_1(\tau)$ is given by the generating series
$$\hat\phi_1(\tau) := \sum_{t\in \Z} \Zh(t,v)\,q^t.$$
\endproclaim

It is conjectured in \cite{\kryII} that the constant $\bold c$ occurring in
the definition of $\Zh(0,v)$ is, in fact, zero.
It may be possible to use recent work of Bruinier and K\"uhn, \cite{\brkuehn}, on the
heights of curves on Hilbert  modular surfaces to show that
that $\langle \hat\o,\hat\o\rangle$ has the predicted value and hence verify this
conjecture.

Some justification for the terminology `arithmetic theta function' is given by the
following result, which is closely related to earlier work of Zagier, \cite{\zagier},
and recent results of Borcherds, \cite{\borchduke}, cf. also \cite{\mcgraw}.

\proclaim{Theorem 1}
The arithmetic theta function $\hat\phi_1(\tau)$ is a
(nonholomorphic) modular form of weight $\frac32$,
valued in $\CH^1(\Cal M)$,
for a subgroup $\Gamma'\subset \SL_2(\Z)$.
\endproclaim

The proof of Theorem 1 depends on Borcherd's result \cite{\borchduke} and on the modularity of various complex
valued $q$-expansions obtained by taking height pairings of $\hat\phi_1(\tau)$ with other classes in $\CH^1(\Cal
M)$. We now describe some of these in terms of values and derivatives of a certain Eisenstein series,
\cite{\kryII}, of weight $\frac32$
$$\Cal E_1(\tau,s,D(B)) = \sum_{\gamma\in \Gamma_\infty\back \SL_2(\Z)}
(c\tau+d)^{-\frac32}|c\tau+d|^{-(s-\frac12)} \,v^{\frac12(s-\frac12)}\,\P_1(s,\gamma,D(B)),$$ associated to $B$
and the lattice $L$, and normalized so that it is invariant under $s\mapsto -s$. The main result of joint work
with M. Rapoport and T. Yang is the following: \proclaim{Theorem 2} {\rm (\cite{\kryII})} {\rm(}i\,{\rm)}
$$\Cal E_1(\tau,\frac12;D(B)) = \deg(\hat\phi_1(\tau))
= \sum_{t} \deg_\Q(\Zh(t,v))\,q^t.$$
{\rm(}ii\,{\rm)}
$$\Cal E'_1(\tau,\frac12;D(B)) = \langle\,
\hat\phi_1(\tau),\hat\o\,\rangle = \sum_{t} \langle\, \Zh(t,v),\hat\o\,\rangle\,q^t.$$
\endproclaim
Note that this result expresses the Fourier coefficients of the first two terms in the
Laurent expansion at the point $s=\frac12$ of
the Eisenstein series $\Cal E_1(\tau,s;D(B))$
in terms of the geometry and the arithmetic geometry of cycles on $\Cal M$.

Next consider the image of
$$\hat\phi_1(\tau) - \Cal E_1(\tau,\frac12;D(B))\cdot\deg(\hat\o)^{-1}\cdot\hat\o$$
in $CH^1(\Cal M_\Q)$, the usual Chow group of the generic fiber.
By (i) of Theorem 2,
it lies in the Mordell-Weil space $CH^1(\Cal M_\Q)^0\otimes\C
\simeq \text{\rm Jac}(M)(\Q)\otimes_\Z\C$.
In fact, it is essentially the generating function
defined by Borcherds, \cite{\borchduke},
for the Shimura curve $M$, and hence is
a holomorphic modular of weight $\frac32$. For the case of modular curves, such a
modular generating
function, whose coefficients are Heegner points, was introduced by Zagier, \cite{\zagier}.
By the Hodge index theorem for $\CH^1(\Cal M)$, \cite{\bost}, the proof of Theorem 1
is completed by showing that the pairing of $\hat\phi_1(\tau)$ with
each class of the form $(Y_p,0)$, for $Y_p$ a component of the fiber $\Cal M_p$,
$p\mid D(B)$ and each class of the form $(0,\phi)$, where $\phi\in C^\infty(M(\C))$, is
modular.

\specialhead \noindent \boldlarge 2.2.  0-cycles\endspecialhead We next consider a generating function for
$0$-cycles  on $\Cal M$. Recall that the arithmetic Chow group $\CH^2(\Cal M)$,
 with real coefficients,
is generated by pairs $(\Cal Z,g)$, where $\Cal Z$ is a real linear combination of $0$-cycles on $\Cal M$ and $g$
is a real smooth $(1,1)$-form on $\Cal M(\C)$. In fact, the arithmetic degree map, as defined in \cite{\bost},
$$\degh:\CH^2(\Cal M)\rightarrow \R, \qquad
\degh((\Cal Z,g)) = \sum_i n_i\log|k(P_i)| + \frac12\int_{\Cal M(\C)} g,$$
where $\Cal Z = \sum_i n_i\,P_i$ for closed points $P_i$ of $\Cal M$
with residue field $k(P_i)$, is an isomorphism.

Let $\tau = u+iv\in \H_2$, the Siegel space of genus $2$, and for $T\in \Sym_2(\Z)$,
let $q^T = e^{2\pi i \tr(T\tau)}$. To define the generating series
$$\hat\phi_2(\tau) = \sum_{T\in \Sym_2(\Z)} \Zh(T,v)\,q^T,$$
we want to define classes $\Zh(T,v)\in \CH^2(\Cal M)$
for each $T\in \Sym_2(\Z)$ and $v\in \Sym_2(\R)_{>0}$.

We begin by considering cycles on $\Cal M$ which are defined by imposing pairs of endomorphisms. For $T\in
\Sym_2(\Z)_{>0}$ a positive definite integral symmetric matrix, let  $\Cal Z(T)$ be the moduli stack over $\Cal M$
consisting of triples $(A,\iota,\bold x)$ where $(A,\iota)$ is as before, and $\bold x= [x_1,x_2]\in V(A,\iota)^2$
is a pair of special endomorphisms with matrix of inner products $Q(\bold x) = \frac12 ((x_i,x_j)) = T$. We call
$T$ the fundamental matrix of the triple $(A,\iota,\bold x)$. The following result of joint work with M. Rapoport
describes the cases in which $\Cal Z(T)$ is, in fact, a $0$-cycle on $\Cal M$. \proclaim{Proposition 2} {\rm
(\cite{\krinvent})} Suppose that $T\in \Sym_2(\Z)_{>0}$. {\rm(}i\,{\rm)} The cycle $\Cal Z(T)$ is either empty or
is supported in the set of supersingular points in a fiber $\Cal M_p$ for a  unique prime $p$ determined by $T$.
In particular, $\Cal Z(T)_\Q=\emptyset$. The prime $p$ is determined by the condition that $T$ is represented by
the ternary quadratic space $V^{(p)} =\{\ x\in B^{(p)}\mid \tr(x)=0\ \}$, with $Q^{(p)}(x) = -x^2,$ where
$B^{(p)}$ is the definite quaternion algebra over $\Q$ with $B^{(p)}_\ell\simeq B_\ell$ for all primes $\ell\ne
p$. If there is no such prime, then $\Cal Z(T)$ is empty. \hfb {\rm(}ii\,{\rm)} {\rm (\, $T$ regular)} Let $p$ be
as in {\rm(}i\,{\rm)}. Then, if $p\nmid D(B)$ or if $p\mid D(B)$ but $p^2\nmid T$, then $\Cal Z(T)$ is a $0$-cycle
in $\Cal M_p$. \hfb {\rm(}iii\,{\rm)} {\rm (\, $T$ irregular)} Let $p$ be as in {\rm(}i\,{\rm)}. If $p\mid D(B)$
and $p^2\mid T$, then $\Cal Z(T)$ is a union, with multiplicities, of components of $\Cal M_p$, cf.
\cite{\krinvent}, 176.
\endproclaim

For $T\in \Sym_2(\Z)_{>0}$ regular, as in (ii) of Proposition 2, we let
$$\Zh(T,v) := \Zh(T) = (\Cal Z(T),0)\in \CH^2(\Cal M).$$

For $T = \pmatrix t_1&m\\m&t_2\endpmatrix\in \Sym_2(\Z)_{>0}$ irregular, we use the
results of \cite{\krinvent}, section 8 (where the quadratic
form is taken with the opposite sign).
We must therefore assume that $p\ne2$, although the results of the
 appendix to section 11 of
\cite{\kryII} suggest that it should be possible to eliminate this restriction.
In this case, the vertical cycle $\Cal Z(T)$ in the fiber $\Cal M_p$
is the union of those connected components
of the
intersection $\Cal Z(t_1)\times_{\Cal M}\Cal Z(t_2)$ where the `fundamental matrix',
\cite{\krinvent}, is equal to $T$. Here
$\Cal Z(t_1)$ and $\Cal Z(t_2)$ are the codimension $1$ cycles defined earlier.
Note that,
by Proposition 1, they can share some vertical components. We base
change to $\Z_p$ and set
$$\Zh(T,v) := \chi(\Cal Z(T),\Cal O_{\Cal Z(t_1)}\otimes^{\Bbb L}
\Cal O_{\Cal Z(t_2)})\cdot \log(p)\quad \in \R\simeq \CH^2(\Cal M),$$ where $\chi$ is the Euler-Poincar\'e
characteristic of the derived tensor product of the structure sheaves $\Cal O_{\Cal Z(t_1)}$ and $\Cal O_{\Cal
Z(t_2)}$, cf. \cite{\krinvent}, section 4. Note that the same definition could have been used in the regular case.

Next we consider nonsingular $T\in \Sym_2(\Z)$ of signature $(1,1)$ or $(0,2)$.
In this case, $\Cal Z(T)$ is empty, since the quadratic form on $V(A,\iota)$
is positive definite, and our
`cycle' should be viewed as `vertical  at infinity'. For a pair of vectors $\bold x =
[x_1,x_2]\in V(\Q)^2$ with  nonsingular matrix of inner products
$Q(\bold x) = \frac12\,\big((x_i,x_j)\big)$, the quantity
$$\L(\bold x) := \int_{D} \xi(x_1)*\xi(x_2),$$
where $\xi(x_1)*\xi(x_2)$ is the $*$-product of the Green functions $\xi(x_1)$ and $\xi(x_2)$, \cite{\gsihes}, is
well defined and depends only on $Q(\bold x)$. In addition, $\L(\bold x)$ has the following remarkable invariance
property.

\proclaim{Theorem 3} {\rm (\cite{\annals, Theorem~11.6})} For $k\in O(2)$, $\L(\bold x\cdot k) =\L(x)$.
\endproclaim

For $T\in \Sym_2(\Z)$ of signature $(1,1)$ or $(0,2)$ and for
$v\in\Sym_2(\R)_{>0}$, choose $a\in GL_2(\R)$ such that $v= a\,{}^ta$, and define
$$\Zh(T,v) := \sum_{\matrix \scr \bold x\in L^2, \  \scr Q(\bold x)=T,
\scr\mod \Gamma\endmatrix} \L(\bold x\,a)\quad \in\R\simeq\CH^2(\Cal M).$$
Here $L=O_B\cap V$ and $\Gamma = O_B^\times$, as before. Note that the invariance
property of Theorem 3 is required to make the right side independent of the choice of
$a$.

We omit the definition of the terms for singular $T$'s, cf. \cite{\kbourb}.

By analogy with Theorem 1, we conjecture that, with this definition,
the generating series $\hat\phi_2(\tau)$
is the $q$-expansion of a Siegel modular form of weight $\frac32$ for a subgroup
$\Gamma'\subset
\roman{Sp}_2(\Z)$.  More precisely, there is a normalized Siegel Eisenstein series $\Cal
E_2(\tau,s;D(B))$ of weight $\frac32$ attached to $B$, \cite{\annals}.
\proclaim{Conjecture 1}
$$\Cal E'_2(\tau,0;D(B)) \overset{?}\to{=}\ \hat\phi_2(\tau).\tag $C1$ $$
This amounts to the family of identities
$$\Cal E'_{2,T}(\tau,0;D(B)) \overset{?}\to{=}\ \Zh(T,v)\,q^T\tag $C1_T$ $$
on Fourier coefficients, for all $T\in \Sym_2(\Z)$.
Here the isomorphism $\degh$\! is being used.
\endproclaim

\proclaim{Theorem 4} {\rm (\cite{\annals}, \cite{\krinvent})} The Fourier coefficient identity {\rm($C1_T$)} holds
in the following cases:\hfb {\rm(}i\,{\rm)} $T\in \Sym_2(\Z)$ is not represented by $V$ or by any of the spaces
$V^{(p)}$  of Proposition 2. \hfb (In this case both $\Zh(T,v)$ and $\Cal E'_{2,T}(\tau,0;D(B))$ are zero.) \hfb
{\rm(}ii\,{\rm)} $T\in \Sym_2(\Z)_{>0}$ is {\it regular} and $p\nmid 2D(B)$, \cite{\annals}. \hfb
{\rm(}iii\,{\rm)} $T\in \Sym_2(\Z)_{>0}$ is irregular with $p\ne2$, or {\it regular} with $p\mid D(B)$ and
$p\ne2$, \cite{\krinvent}. \hfb {\rm(}iv\,\rm{)} $T\in \Sym_2(\Z)$ is nonsingular of signature $(1,1)$ or $(0,2)$,
\cite{\annals}.
\endproclaim

 Theorem 4 is proved by a direct computation of both
sides of ($C1_T$). In case (ii), the computation of the Fourier coefficient $\Cal E'_{2,T}(\tau,0;D(B))$ depends
on the formula of Kitaoka, \cite{\kitaoka}, for the local representation densities $\a_p(S,T)$ for the given $T$
and a variable unimodular $S$. The computation of $\Zh(T,v) = \degh\!((\Cal Z(T),0))$ depends on a special case of
a result of Gross and Keating, \cite{\grosskeating}, about the deformations of a triple of isogenies between a
pair of $p$-divisible formal groups of dimension $1$ and height $2$ over $\bar\Bbb F_p$. Their result is also
valid for $p=2$, so it should be possible to extend (ii) to the case $p=2$ by extending the result of Kitaoka.

 In case (iii), an explicit formula for the quantity $\chi(\Cal Z(T),\Cal O_{\Cal Z(t_1)}\otimes^{\Bbb L} \Cal O_{\Cal
Z(t_2)})$ is obtained in \cite{\krinvent} using $p$-adic uniformization. The analogue of Kitaoka's result is a
determination of $\a_p(S,T)$ for arbitrary $S$ due to T. Yang, \cite{\yangden}. In both of these results, the case
$p=2$ remains to be done.

 Case (iv) is proved by directly relating the function $\Lambda$, defined via the
$*$-product to the derivative at $s=0$ of the confluent hypergeometric function of a matrix argument defined by
Shimura, \cite{\shimura}. The invariance property of Theorem 3 plays an essential role. The case of signature
$(1,1)$ is done in \cite{\annals}; the argument for signature $(0,2)$ is the same.

A more detailed sketch of the proofs can be found in \cite{\kbourb}.

 As part of ongoing joint work with
M. Rapoport and T. Yang, 
the verification of ($C1_T$) for singular $T$ of rank $1$
is nearly complete.

\specialhead \noindent \boldLARGE 3. Higher dimensional
examples \endspecialhead

So far, we have discussed the generating functions $\hat\phi_1(\tau)\in \CH^1(\Cal M)$ and
$\hat\phi_2(\tau)\in \CH^2(\Cal M)$ attached to the arithmetic surface $\Cal M$,
and the connections of these series to derivatives of Eisenstein series.
There should be analogous series defined as generating functions for arithmetic cycles
for the Shimura varieties attached to rational quadratic spaces $(V,Q)$ of signature
$(n,2)$. At present there are several additional examples, all based on
the accidental isomorphisms for small values of $n$, which allow us to
identify the Shimura varieties in question with moduli spaces of abelian
varieties with specified polarization and endomorphisms. Here we briefly sketch
what one hopes to obtain and indicate what is known so far. The results
here are joint work with M. Rapoport.

{\bf Hilbert-Blumenthal varieties ($n=2$), \cite{\krHB}.} When the rational quadratic space $(V,Q)$ has signature
$(2,2)$, the associated Shimura variety $M$ is a quasi-projective surface with a canonical model over $\Q$. There
is a model $\Cal M$ of $M$ over $\Spec(\Z[N^{-1}])$ defined as the moduli scheme for collections
$(A,\l,\iota,\bar\eta)$ where $A$ is an abelian scheme of relative dimension $8$ dimension with polarization $\l$,
level structure $\bar\eta$, and an action of $O_C\otimes\OK$, where $O_C$ is a maximal order in the Clifford
algebra $C(V)$ of $V$ and $\OK$ is the ring of integers in the quadratic field $\kay = \Q(\sqrt{d})$ for $d=
\roman{discr(V)}$, the discriminant field of $V$, \cite{\krHB}. Again, a space  $V(A,\iota) =
V(A,\lambda,\iota,\bar\eta)$ of special endomorphisms is defined; it is a $\Z$-module of finite rank equipped with
a positive definite quadratic form $Q$. For $T\in \Sym_r(\Z)$, we let $\Cal Z(T)$ be the locus of
$(A,\l,\iota,\bar\eta,\bold x)$'s where $\bold x = [x_1,\dots,x_r]$, $x_i\in V(A,\iota)$  is a collection of $r$
special endomorphisms with matrix of inner products $Q(\bold x) = \frac12\big((x_i,x_j)\big) =T$.

One would {\it like} to define a family of generating functions according
to the following conjectural chart. Again there is a metrized Hodge line
bundle $\hat\o\in \CH^1(\Cal M)$.
\smallskip
\settabs 8\columns \+$r=1,$&\hskip -10pt$\Cal Z(t)_\Q =\roman{HZ-curve,}$&&\hskip -10pt $\hat\phi_1(\tau) = \hat\o
+ ? +\sum_{t\ne0} \Zh(t,v)\,q^t,$&&& $\langle \hat\phi_1(\tau),\hat\o^2\rangle\overset{?}\to{ =}\ \Cal
E'_1(\tau,1).$ \cr \+$r=2,$&\hskip -10pt$\Cal Z(t)_\Q = $0$-\roman{cycle,}$&&\hskip -10pt $\hat\phi_2(\tau) =
\hat\o^2 + ? +\sum_{T\ne0} \Zh(T,v)\,q^T,$&&& $\langle \hat\phi_2(\tau),\hat\o\rangle\overset{?}\to{ =}\ \Cal
E'_2(\tau,\frac12).$ \cr \+$r=3,$&\hskip -10pt$\Cal Z(T)_\Q =\emptyset,$&&\hskip -10pt $\hat\phi_3(\tau) =
\hat\o^3 + ? +\sum_{T\ne0} \Zh(T,v)\,q^T,$&&& $\degh\!\hat\phi_3(\tau)\overset{?}\to{ =}\ \Cal E'_3(\tau,0).$ \cr
Here, the generating function $\hat\phi_r(\tau)$ is valued in $\CH^r(\Cal M)$, the $r$th arithmetic Chow group,
$\Cal E_r(\tau,s)$  is a certain normalized Siegel Eisenstein series of genus $r$, and the critical value of $s$
in the identity in the last column is the Siegel-Weil point $s_0 = \frac12(\dim(V) - r-1)$. Of course, one would
like the $\hat\phi_r(\tau)$'s to be Siegel modular forms of genus $r$ and weight $2$.

There are many technical problems which must be overcome to obtain such results. For example, one would like to
work with a model over $\Spec(\Z)$. If $V$ is anisotropic, then $M$ is projective, but if $V$ is isotropic, e.g.,
for the classical Hilbert-Blumenthal surfaces where it has $\Q$-rank $1$, then one must compactify. Since the
metric on $\hat\o$ is singular at the boundary a more general version of the Gillet-Soul\'e theory, currently
being developed by Burgos, Kramer and K\"uhn, \cite{\bkk}, \cite{\kuehn}, will be needed.

Nonetheless, the chart suggests many identities
which can in fact be checked rigorously. For example, there are again
rational quadratic spaces $V^{(p)}$ of dimension $4$ and signature $(4,0)$ obtained
by switching the Hasse invariant of $V$ at $p$.
\proclaim{Theorem 5} {\rm \cite{\krHB}, \cite{\kbourb}.} {\rm(}i\,{\rm)} If $T\in
\Sym_3(\Z)_{>0}$  is not represented by any of the $V^{(p)}$'s, then $\Cal Z(T)=\emptyset$
and $\Cal E'_{3,T}(\tau,0)=0$. \hfb
{\rm(}ii\,{\rm)} If $T\in \Sym_3(\Z)_{>0}$ is represented by $V^{(p)}$ where $p$ is
a prime of good reduction
split in $\kay$, then $\Cal Z(T)$ is a $0$-cycle in $\Cal M_p$ and
$$\degh\!((\Cal Z(T),0))\, q^T = \Cal E'_{3,T}(\tau,0).\tag $\star$ $$
{\rm(}iii\,{\rm)} If $T\in \Sym_3(\Z)_{>0}$ is represented by $V^{(p)}$ and $p$ is
a prime of good reduction inert in $\kay$,
then $\Cal Z(T)$ is a $0$-cycle in $\Cal M_p$
if and only if
$p\nmid T$. If this is the case, then the Fourier coefficient identity ($\star$)
again holds. If $p\mid T$, then $\Cal Z(T)$ is a union of components of the
supersingular locus of $\Cal M_p$.
\endproclaim

Finally, say if $V$ is anisotropic, one can consider the image
$\roman{cl}(\hat\phi_r(\tau))\in H^{2r}(M,\C)$ of
$\hat\phi_r(\tau)$ in the usual (Betti) cohomology of $\Cal M(\C)$.
Of course, $\roman{cl}(\hat\phi_3(\tau))=0$ for degree reasons.
Joint work with J. Millson on generating functions for cohomology classes
of special cycles yields:
\proclaim{Theorem 6} {\rm (\cite{\kmihes}, \cite{\duke}, \cite{\kbourb})} Suppose that $V$ is anisotropic.
{\rm(}i\,{\rm)} $\roman{cl}(\hat\phi_r(\tau))$ is a Siegel modular form of genus $r$ and weight $2$ valued in
$H^{2r}(M,\C)$. \hfb {\rm(}ii\,{\rm)} For the cup product pairing, $\big(\ \roman{cl}(\hat\phi_r(\tau)),
\roman{cl}(\hat\o)\ \big) = \Cal E_r(\tau,s_0)$, where  $s_0 = \frac12(3-r)$.
\endproclaim
Part (ii) here generalizes (i) of Theorem 2 above, so that, again, the
value at $s_0$ of the Eisenstein series $\Cal E_r(\tau,s)$ involves
the complex geometry, while,
conjecturally,  the second term involves the height pairing.

{\bf Siegel modular varieties ($n=3$), \cite{16}.}
Here, an integral model $\Cal M$ of the Shimura variety $M$ attached to a rational
quadratic space of signature $(3,2)$ can be obtained as a moduli space
of polarized abelian varieties of dimension $16$ with an action of
a maximal order $O_C$ in the Clifford algebra of $V$. We just give the
relevant {\it conjectural} chart:
\smallskip
\settabs 8\columns \+$r=1,$&\hskip -10pt$\Cal Z(t)_\Q ={\textstyle\text{Humbert}\atop\textstyle \lower
2pt\hbox{\text{surface}}} $,&&\hskip -10pt $\hat\phi_1(\tau) = \hat\o + ? +\sum_{t\ne0} \Zh(t,v)\,q^t,$&&&
$\langle \hat\phi_1(\tau),\hat\o^3\rangle\overset{?}\to{ =}\ \Cal E'_1(\tau,\frac32).$ \cr \+$r=2,$&\hskip
-10pt$\Cal Z(t)_\Q = \roman{curve}$&&\hskip -10pt $\hat\phi_2(\tau) = \hat\o^2 + ? +\sum_{T\ne0}
\Zh(T,v)\,q^T,$&&& $\langle \hat\phi_2(\tau),\hat\o^2\rangle\overset{?}\to{ =}\ \Cal E'_2(\tau,1).$ \cr
\+$r=3,$&\hskip -10pt$\Cal Z(T)_\Q =0$-cycle,&&\hskip -10pt $\hat\phi_3(\tau) = \hat\o^3 + ? +\sum_{T\ne0}
\Zh(T,v)\,q^T,$&&& $\langle \hat\phi_2(\tau),\hat\o\rangle\overset{?}\to{ =}\ \Cal E'_3(\tau,\frac12).$ \cr
\+$r=4,$&\hskip -10pt$\Cal Z(T)_\Q =\emptyset,$&&\hskip -10pt $\hat\phi_4(\tau) = \hat\o^4+ ? +\sum_{T\ne0}
\Zh(T,v)\,q^T,$&&& $\degh\!\hat\phi_4(\tau)\overset{?}\to{ =}\ \Cal E'_4(\tau,0).$ \cr Here the Eisenstein series
and, conjecturally, the generating functions $\hat\phi_r(\tau)$ have weight $\frac52$, and the values of the
Eisenstein series should be related to the series $\roman{cl}(\hat\phi_r(\tau))$. In the case of a prime $p$ of
good reduction a model of $M$ over $\Spec(\Z_p)$ is defined in \cite{\krsiegel}, and cycles are defined by
imposing special endomorphisms. For example, for $r=4$, the main results of \cite{\krsiegel} give a criterion for
$\Cal Z(T)$ to be a $0$-cycle in a fiber $\Cal M_p$ and show that, when this is the case, then $\degh\!((\Cal
Z(T),0))\,q^T = \Cal E'_{4,T}(\tau,0)$. The calculation of the left hand side is again based on the result of
Gross and Keating mentioned in the description of the proof of Theorem 4 above. This provides some evidence for
the last of the derivative identities in the chart.

\redefine\vol{\oldvol}

\specialhead \noindent \boldLARGE References \endspecialhead

\widestnumber\key{AA}

\parskip=0pt

\ref\key{\borchduke}
\by R. Borcherds
\paper\rm  The Gross-Kohnen-Zagier theorem in higher dimensions
\jour\it Duke Math. J.
\yr 1999
\vol 97
\pages 219--233
\endref

\ref\key{\bost}
\by J.-B. Bost
\paper\rm  Potential theory and Lefschetz theorems for arithmetic surfaces
\jour\it Ann. Sci. \'Ecole Norm. Sup.
\yr 1999
\vol 32
\pages 241--312
\endref

\ref\key{\boutotcarayol}
\by J.-F. Boutot and H. Carayol
\paper\rm  Uniformisation $p$-adique des courbes de Shimura
\inbook in Ast\'erisque, vol. {\bf196--197}
\yr 1991
\pages 45--158
\endref

\ref\key{\brkuehn}
\by J. H. Bruinier, J. Burgos, and U. K\"uhn
\paper\rm  in preparation
\endref

\ref\key{\bkk}
\by J. Burgos, J. Kramer and U. K\"uhn
\paper\rm  in preparation
\jour\it \yr
\vol
\pages
\endref

\ref\key{\gsihes}
\by H. Gillet and C. Soul\'e
\paper\rm  Arithmetic intersection theory
\jour\it Publ. Math. IHES
\yr 1990
\vol 72
\pages 93--174
\endref

\ref\key{\grosskeating}
\by B. H. Gross and K. Keating
\paper\rm  On the intersection of modular correspondences
\jour\it Invent.\ Math.
\vol 112
\yr 1993
\pages 225--245
\endref

\ref\key{\kitaoka}
\by Y. Kitaoka
\paper\rm  A note on local densities of quadratic forms
\jour\it Nagoya Math. J.
\vol 92
\yr 1983
\pages 145--152
\endref

\ref\key{\duke}
\by S. Kudla
\paper\rm  Algebraic cycles on Shimura varieties of orthogonal type
\jour\it Duke Math. J.
\yr 1997
\vol 86
\pages 39--78
\endref

\ref\key{\annals}
\bysame
\paper\rm  Central derivatives of Eisenstein series and height pairings
\jour\it Ann. of Math.
\vol 146
\yr 1997
\pages 545-646
\endref

\ref\key{\kbourb}
\bysame
\paper\rm  Derivatives of Eisenstein series and generating functions for arithmetic cycles
\inbook S\'em. Bourbaki n${}^o$ 876
\bookinfo Ast\'erisque
\vol {\bf 276}
\yr 2002
\pages 341--368
\endref

\ref\key{\kudlamsri}
\bysame
\paper\rm  Special cycles and derivatives of Eisenstein series
\jour\it Proc. of MSRI Workshop on Heegner points (to appear)
\endref

\ref\key{\kmihes}
\by S. Kudla and J. Millson
\paper\rm  Intersection numbers of cycles on locally symmetric spaces and Fourier
coefficients of holomorphic modular forms in several complex variables
\jour\it Publ. Math. IHES
\vol 71
\yr 1990
\pages 121--172
\endref

\ref\key{\krHB}
\by S. Kudla and M. Rapoport
\paper\rm  Arithmetic Hirzebruch--Zagier cycles
\jour\it J. reine angew. Math.
\vol 515
\yr 1999
\pages 155--244
\endref

\ref\key{\krinvent}
\bysame
\paper\rm  Height pairings on Shimura curves and $p$-adic unformization
\jour\it Invent. math.
\yr 2000
\vol 142
\pages 153--223
\endref

\ref\key{\krsiegel}
\bysame
\paper\rm  Cycles on Siegel threefolds and derivatives of Eisenstein series
\jour\it Ann. Scient. \'Ec. Norm. Sup.
\vol 33
\yr 2000
\pages 695--756
\endref

\ref\key{\tiny}
\by S. Kudla, M. Rapoport and T. Yang
\paper\rm  On the derivative of an Eisenstein series of weight 1
\jour\it Int. Math. Res. Notices, No.7
\yr 1999
\pages 347--385
\endref

\ref\key{\kryII}
\bysame
\paper\rm  Derivatives of Eisenstein series and Faltings heights
\jour\it preprint
\yr 2001
\vol
\pages
\endref

\ref\key{\kuehn}
\by U. K\"uhn
\paper\rm  Generalized arithmetic intersection numbers
\jour\it J. reine angew. Math.
\yr 2001
\vol 534
\pages 209--236
\endref

\ref\key{\mcgraw}
\by W. J. McGraw
\paper\rm  On the rationality of vector-valued modular forms
\jour\it preprint
\yr 2001
\endref

\ref\key{\shimura}
\by G. Shimura
\paper\rm  Confluent hypergeometric functions on tube domains
\jour\it Math. Annalen
\vol 260
\yr 1982
\pages 269--302
\endref

\ref\key{\yangden}
\by T. Yang
\paper\rm  An explicit formula for local densities of quadratic forms
\jour\it J. Number Theory
\vol 72
\yr 1998
\pages 309--356
\endref

\ref\key{\yangiccm}
\bysame
\paper\rm  The second term of an Eisenstein series
\jour\it Proc. of the ICCM, (to appear)
\endref

\ref\key{\yangmsri}
\bysame
\paper\rm  Faltings heights and the derivative of Zagier's Eisenstein series
\jour\it Proc. of MSRI workshop on Heegner points, preprint (2002)
\endref

\ref\key{\zagier} \by D. Zagier \paper\rm  Modular points, modular curves, modular surfaces and modular forms
\inbook Lecture Notes in Math.~1111 \yr 1985, 225--248 \publ Springer \publaddr Berlin
\endref

\enddocument